\input amstex
\documentstyle{amsppt}
\input bull-ppt
\keyedby{bull298e/pah}

\topmatter           
\cvol{27}
\cvolyear{1992}
\cmonth{July}
\cyear{1992}
\cvolno{1}
\cpgs{154-158}
\title  Generalizing The Hyperbolic Collar Lemma  
\endtitle
\author Ara Basmajian
\endauthor
\address Mathematics department, University of Oklahoma, 
Norman, Oklahoma 73019\endaddress
\ml abasmajian\@nsfuvax.math.uoknor.edu\endml
\date October 25, 1991 and, in revised form, January 28, 
1992\enddate
\abstract We discuss two generalizations of the collar 
lemma.
 The first is the stable neighborhood theorem
which says that a (not necessarily simple) 
closed geodesic in a hyperbolic surface
has a \lq\lq stable neighborhood\rq\rq  whose width only 
depends
on the length of the geodesic. As an application, we show 
that there is
a lower bound for the length of a closed geodesic having 
crossing number $k$ on a hyperbolic surface. This lower
bound tends to infinity with $k$.

Our second generalization is to totally geodesic 
hypersurfaces of
hyperbolic manifolds. Namely, we construct a tubular 
neighborhood
function and show that an embedded closed 
totally geodesic hypersurface in a hyperbolic manifold has 
a 
tubular neighborhood whose width only depends on the area of
the hypersurface (and hence not on the geometry of the 
ambient manifold).
The implications of this result for volumes of hyperbolic 
manifolds is
discussed. We also derive a (hyperbolic) quantitative 
version
of the Klein-Maskit combination theorem (in all dimensions) 
for free products of fuchsian groups.
Using this last theorem, we construct examples  to 
illustrate the qualitative
sharpness of the tubular neighborhood function.\endabstract
\endtopmatter

\document
\heading 1. Introduction\endheading
It is natural to look for universal properties of discrete 
subgroups
of a fixed lie group. The lie group we are concerned with 
is the
orientation preserving isometries of real hyperbolic space
$\text{Isom}^+(\Bbb H^n)$. A fundamental universal 
property in dimension
three is the Jorgensen inequality [Jo].  
Another celebrated property is the 
Margulis lemma [Th] (actually, this lemma holds for a 
bigger class of lie groups). We would like to focus on the
well-known collar lemma due to Keen [Ke]; the sharp form is 
due to Buser [Bu]. There are many different versions of 
this 
lemma in the literature (cf.\ [Ba1, Be, Ha, M1, Ma, R]).

Let $r:\Bbb R^+ \rightarrow \Bbb R^+$ be the function
$r(x)=\log \coth ({x}/{2})$.
Observe that this function monotonically decreases to zero 
and
satisfies $r^2=1$.
The collar lemma says that any simple closed geodesic of 
length
$\ell$ on a hyperbolic surface (a complete surface of 
constant
 curvature negative one) has a collar neighborhood of width
$r(\frac{\ell}{2})$. Furthermore, a collection of disjoint 
simple
closed geodesics have disjoint collars. This result is 
universal in
the sense that the collar width does not depend on the 
underlying
hyperbolic structure of the surface, only on the length of 
the
geodesic. 

We would like to discuss two generalizations of the collar 
lemma.
In the first generalization (the stable neighborhood 
theorem), we 
remove the hypothesis in the collar lemma 
that the closed geodesic be simple. The concept
of a collar is then replaced by the notion of a stable 
neighborhood.
The second generalization  of the collar lemma 
involves totally geodesic hypersurfaces in hyperbolic 
manifolds 
(complete riemannian manifolds of constant curvature 
negative one).
Simple closed geodesics are replaced by totally geodesic 
embedded
hypersurfaces and the notion of a collar is then replaced 
by a
tubular neighborhood.

\heading 2. The stable neighborhood theorem  \endheading

The \it neighborhood of width $d$, \rm $U_d(\omega)$, 
about a geodesic 
$\omega$ in a hyperbolic surface
is the set of all points  within a distance $d$ from the 
geodesic. 
Suppose $\omega$ is a closed geodesic. Then the neighborhood
$U_d(\omega)$ is said to be \it stable \rm if 
for any two connected (smooth) lifts of $\omega$, say 
$\omega_1$ and $\omega_2$, we have
$$
\omega_1 \cap \omega_2 \neq \emptyset \quad \text{if and 
only if}\quad
U_d(\omega_1)\cap U_d(\omega_2)\neq \emptyset.
$$

\proclaim{Theorem 1 {\rm (The Stable Neighborhood Theorem)}}
 A closed geodesic of length $\ell$ on a hyperbolic
surface has a stable neighborhood of width 
$r(\frac{\ell}{2}).$
Furthermore, two disjoint closed geodesics $\omega_1$ and 
$\omega_2$
on the surface having lengths
$\ell_1$ and $\ell_2$  have disjoint stable neighborhoods 
of widths 
$r(\frac{\ell_1}{2})$
and $r(\frac{\ell_2}{2})$ respectively, if $\omega_1$ and 
$\omega_2$
 are separated by a disjoint union
of simple closed geodesics. 
\endproclaim

The above separation condition  is necessary,  since there 
are examples of 
closed geodesics that get arbitrarily close to any fixed 
boundary
geodesic on a pair of pants. 

Hempel, Nakanishi, and Yamada [He, N, Y1, Y2]) have 
independently
shown that there is a universal lower bound for the length 
of a 
nonsimple closed geodesic on a hyperbolic surface. As a 
consequence
of the stable neighborhood theorem we have

\proclaim{Corollary 2} 
There exists an increasing
sequence $M_k$ \RM(for $k=1,2,3,\dots)$, tending to infinity
so that if $\omega$ is a closed geodesic with 
self-intersection 
number $k$, then $\ell (\omega)>M_k$.
Thus the length of a closed geodesic gets arbitrarily large
as its self-intersection gets large.  
\endproclaim

\heading 3. The tubular neighborhood theorem \endheading

In [Ba2], it was shown that associated to any
totally geodesic hypersurface in a hyperbolic manifold 
there exists
a spectrum of numbers called the orthogonal spectrum. 
This spectrum is essentially the lengths of orthogonals 
that start
and end in the hypersurface. The
following can be thought of as a geometric study of the 
first element 
in that spectrum.

Let $V_n(r)$ be the volume of the $n$-dimensional 
hyperbolic ball of radius
$r$. The function
$$
c_n(A)=\frac{1}{2}(V_{n-1}\circ r )^{-1} (A)
$$ is called the \it $n$-dimensional tubular neighborhood 
function. \rm
Observe that this function is monotone decreasing and 
tends to zero as
$A$ goes to infinity. The \it area \rm of a hypersurface 
is with
respect to the induced metric from the ambient hyperbolic 
manifold.

The following theorem shows that a closed 
embedded totally geodesic hypersurface in a hyperbolic 
manifold
has a tubular neighborhood whose width  only depends on 
the area 
of the hypersurface.

\proclaim{Theorem 3 {\rm (The Tubular Neighborhood 
Theorem)}}
Suppose $M^n$ is a hyperbolic manifold 
 containing
$\Sigma$, an embedded closed totally geodesic hypersurface 
of area $A$. 
Then $\Sigma$
has a tubular neighborhood of width $c_n(A).$ That is, the 
set of points
$$
\{x \in M : \roman d(x,\Sigma) < c_n (A)\}
$$
is isometric to the product $\Sigma \times 
(-c_n(A),c_n(A)).$

Furthermore, any disjoint set of such hypersurfaces has 
disjoint
tubular neighborhoods.
\endproclaim

The main idea in the proof of the tubular neighborhood 
theorem is that
the hypersurface in question must contain an embedded disc 
of radius
$r(d)$, where $d$ is the length of the shortest 
common orthogonal from $\Sigma$ to
itself.

It is well known that there exists a constant $a$, so that 
if
$M$ is a hyperbolic three manifold containing $K$
rank two cusps then the volume of $M$ is bigger than $Ka$ 
(see [Th]).
Analyzing volumes of tubular neighborhoods, we can 
prove a version of this result for totally geodesic 
hypersurfaces in hyperbolic manifolds of all dimensions.
Specifically, we have

\proclaim{Theorem 4} There exist positive constants $a_n$, 
for
each $n \geq 3$, depending only on the dimension $n$, so 
that if $M^n$ is a 
hyperbolic manifold containing $K$  closed embedded 
disjoint totally 
geodesic hypersurfaces
then,
$$
\roman{Vol} (M) > K a_n.
$$
For hyperbolic three manifolds, $a_3$ can be taken to be 
$\pi\left(\log 2 +\frac{\sqrt 2}{2}\right)$, which is 
approximately
{\rm 4.4.}

In dimensions $n \geq 4$, if $M_i^n$ is a sequence of 
\RM(not necessarily distinct\RM) hyperbolic $n$-manifolds
each containing an embedded totally geodesic closed
hypersurface of area $A_i$,
 where $A_i \rightarrow \infty$, then the volumes 
of the tubular neighborhoods of width $c_n(A_i)$ tend to 
infinity.
In particular, the volumes $\roman{Vol}(M_i) \rightarrow 
\infty$.

In dimension three, there exist examples of totally 
geodesic surfaces in 
hyperbolic three manifolds whose areas get arbitrarily 
large, but whose 
best possible tubular neighborhoods have bounded volume.
\endproclaim

The main lemma needed to prove the lower volume bound in 
Theorem 4 is the rate of growth lemma.

\proclaim{Lemma {\rm (Rate of growth lemma)}} $\Cal V_n$, 
the volume of a
\RM(one-sided\RM) tubular neighborhood of width $c_n(x)$
about a hypersurface having area $x$,
has the following behavior at
infinity,
$$
\lim_{x\to \infty}\Cal V_n(x)=\cases \infty &\text{for 
$n\geq 4$},\\
                                   \pi &\text{for 
$n=3$},\\  
                                    0  &\text{for $n=2$}.
                            \endcases                  
$$

$\Cal V_3$ is monotone increasing and $\Cal V_2$ is 
monotone decreasing.

\endproclaim

We remark that Kojima and Miyamoto [KM] have found the
smallest volume hyperbolic three manifolds with totally 
geodesic
boundary. This volume is about 6.45. 

In dimension three, a rank two cusp has a cusp neighborhood
whose boundary is a flat torus having  shortest closed 
geodesic 
of  length one. Furthermore, if the hyperbolic three 
manifold has more than one
cusp,
these regions will be disjoint. In fact, this is the basis 
of 
Meyerhoff's lower
volume bound for a hyperbolic three manifold containing 
rank two cusps 
[Me]. 

The next corollary extends this to our setting:

\proclaim{Corollary 5} The tubular neighborhood 
of width $c_3(A)$ about a totally
geodesic embedded closed surface of area $A$ in a 
hyperbolic three manifold
is disjoint from the tori on the boundary of
a rank two cusp having shortest closed geodesic of length 
one. In particular, if
the three manifold $M$ has $n$ rank two cusps and $k$ 
disjoint totally geodesic
closed embedded surfaces, then
$$
\roman{Vol}(M)>(\sqrt3/4)n+(4.4)k.
$$
\endproclaim   

This corollary follows easily from the next

\proclaim{Lemma} Suppose $M$ is a hyperbolic three 
manifold containing a
\RM(rank one or two\RM) 
cusp and an embedded closed totally geodesic surface 
$\Sigma$.
Let $\Bbb T^2$ be the boundary torus or annulus having its 
shortest 
closed geodesic of length one. Then the distance
$$
\roman d(\Bbb T^2,\Sigma)>\log 2.
$$
\endproclaim

A \it fuchsian subgroup \rm in the group of isometries 
$\text{Isom}^+(\Bbb H^n)$ is a pair $(F,X)$ where $X$ is a 
hyperbolic
hyperplane and $F$ is a 
discrete  subgroup keeping $X$ and the half-spaces it 
bounds in $\Bbb H^n$
 invariant. The \it injectivity radius \rm of $F$ at $x 
\in X$,
denoted $\text{inj(x)}$, is the
largest hyperbolic disc centered at $x$ whose 
$F$-translates are all 
disjoint. If $x$ is an elliptic fixed (orbifold) point, 
then its 
injectivity radius is zero. There is a unique common 
orthogonal 
between any two disjoint (in $\overline{\Bbb H}^n$) 
hyperplanes.

The following is a hyperbolic (quantitative) version of 
the Klein-Maskit combination theorem
for fuchsian subgroups in all dimensions. Its proof makes 
essential use of the combination theorem 
 in all dimensions (see [Ma2, Ap]).

\proclaim{Theorem 6} Suppose $(F_1,X_1)$ and $(F_2,X_2)$ 
are 
fuchsian subgroups of the full isometry group 
$\roman{Isom}^+(\Bbb H^n)$ with
$X_1$ and $X_2$ disjoint in $\Bbb H^n$. Let $x_1 \in X_1$ 
and 
$x_2 \in X_2$ be the endpoints of the unique common 
perpendicular
between the hyperplanes $X_1$ and $X_2$. If the following 
inequality
holds
$$
r(\roman{inj}(x_1))+r(\roman{inj}(x_2))<\roman d(X_1,X_2), 
\tag *
$$
then the group $G=<F_1,F_2>$ is a discrete  group of the 
second kind
that, abstractly, is the free product of $F_1$ and $F_2$. 
Furthermore,
the hyperbolic $(n-1)$-hypersurfaces
 $\Sigma_i=X_i/F_i$ \RM(\,for $i=1,2)$ are totally 
geodesic boundary hypersurfaces for the hyperbolic manifold 
$N=\Bbb H^n/G$ satisfying $\roman d_N(\Sigma_1,\Sigma_2)=
\roman d(X_1,X_2).$
The group $G$ is torsion-free if and only if $F_1$ and 
$F_2$ are 
torsion-free.

\endproclaim

The above theorem with very few exceptions holds when 
there is 
equality in $(*)$. As a consequence of this, we have

\proclaim{Corollary 7} Suppose $\Sigma_1$ and $\Sigma_2$ 
are hyperbolic 
$(n-1)$-manifolds containing embedded balls of radii $R_1$ 
and $R_2$, 
respectively. Then
there exists a hyperbolic $n$-manifold $N$ having totally 
geodesic
boundary hypersurfaces $\Sigma_1$ and $\Sigma_2$ 
satisfying, 
$\roman d_N(\Sigma_1,\Sigma_2)=r(R_1)+r(R_2)$.\qed
\endproclaim

The dimension three examples in Theorem 4 are constructed 
using 
the above corollary.

The proofs of the stable neighborhood theorem, its 
corollary, and
other results on short nonsimple closed geodesics 
are contained in the paper [Ba3]. The proofs of the
theorems on totally geodesic hypersurfaces,
a more general form of Theorem 6, along with additional
consequences of this approach are contained in [Ba4].

\heading Acknowledgements\endheading
 The author would like to thank
Colin Adams, Boris Apanasov, and Bernard Maskit for helpful
conversations.


\Refs
\ra\key{MMM}

\ref \key{Ap} \by Boris Apanasov \book Discrete groups
in space and uniformization problems \publ Kluwer Academic 
\publaddr Dordrecht, Netherlands
\yr1991
\endref

\ref \key{ Ba1} \by Ara Basmajian \paper Constructing 
pairs of pants
\yr 1990 \vol15 \jour Ann. Acad. Sci. Fenn. Ser. A I Math.
\pages 65--74
\endref


\ref \key{ Ba2} \bysame \paper The orthogonal spectrum
of a hyperbolic manifold \jour Amer. J. Math., (to appear, 
1992)
\endref
\ref \key{ Ba3} \bysame\paper The stable neighborhood
theorem and lengths of closed geodesics 
\jour Proc. Amer. Math. Soc., (to appear)
\endref

\ref \key{ Ba4} \bysame
\paper  Tubular neighborhoods of totally geodesic 
hypersurfaces 
in hyperbolic manifolds 
\paperinfo preprint\endref


\ref \key{ Be} \by Lipman Bers \paper An inequality for 
Riemann surfaces
\jour Differential geometry and Complex analysis, H.E. 
Rauch Memorial volume
(Isaac Chavel and Herschel M. Farkas, eds.), 
Springer-Verlag,
New York, 1985, pp. 87--93
\endref

\ref \key{ Bu} \by Peter Buser \paper The collar theorem 
and examples
\yr1978 \vol 25 \jour Manuscripta Math. \pages 349-357
\endref

\ref\key{ Ha}\by N. Halpern
\paper A proof of the collar lemma \jour Bull.
London Math. Soc. \vol 13\yr1981\pages 141--144\endref



\ref 
\key{ He} \by John Hempel \paper Traces, lengths, and 
simplicity of
loops on surfaces \jour Topology Appl. \vol 18
\yr 1984 \pages 153--161
\endref

\ref 
\key{ Jo} \by Troels Jorgensen \paper on discrete groups
of M\"obius transformations
 \jour Amer. J. Math. \vol 98
\yr1976  \pages 739--749 
\endref

\ref \key{ Ke} \by Linda Keen \paper Collars on Riemann 
surfaces
\jour Discontinuous Groups and Riemann Surfaces, Ann. of 
Math. Stud.,
vol. 79, Princeton Univ. Press, Princeton, NJ, 1974, pp. 
263--268
\endref


\ref \key{ KM} \by S. Kojima and Y. Miyamoto
\paper The smallest hyperbolic 3-manifolds with totally 
geodesic
boundary
\yr 1991 \vol34 \jour J. Differential Geometry  \pages 
175--192
\endref

\ref \key{ Ma1} \by Bernard Maskit \paper Comparison of 
hyperbolic and
extremal lengths \jour Ann. Acad. Sci. Fenn. Ser. A I Math.
\vol 10 \yr 1985 \pages 381--386
\endref

\ref \key{ Ma2} \bysame \book Kleinian groups 
\publ Springer-Verlag\publaddr New York
\yr1988
\endref




\ref \key{ Mt} \by Peter Matelski \paper A compactness 
theorem for 
Fuchsian groups of the second kind \jour Duke Math. J. 
\vol 43
\yr 1976 \pages 829--840
\endref

\ref \key{ Me} \by Bob Meyerhoff  \paper A lower bound
for the volume of hyperbolic \RM 3-manifolds \yr 1987 \vol39
\jour Canad. J. Math. \pages 1038--1056 
\endref


\ref \key{ N} \by Toshihiro Nakanishi \paper The lengths 
of the closed
geodesics on a Riemann surface with self-inter\-section 
\jour Tohoku Math. J. (2)
\vol 41 \yr 1989 \pages  527--541
\endref

\ref \key{ R} \by Burton Randol \paper Cylinders in 
Riemann surfaces
\jour Comment. Math. Helv. \vol 54 \yr 1979 \pages 1--5
\endref


\ref \key{ Th} \by William Thurston \book The geometry and
topology of \RM3-manifolds, lecture notes \publ Princeton 
University
\yr 1977
\endref


\ref \key{ Y1} \by Akira Yamada \paper On Marden's universal
constant of Fuchsian groups  \yr 1981 \vol 4 \jour Kodai 
Math. J.
\pages 266--277 
\endref 


\ref \key{ Y2} \by Akira Yamada \paper On Marden's universal
constant of Fuchsian groups {\rm II}\jour J. Analyse 
Math.\vol 41\yr 1982
\pages 234--248
\endref 

\endRefs
\enddocument